\documentclass{article}

\usepackage{amsmath,amsthm,amssymb}
\usepackage{color}
\usepackage{graphicx}

\usepackage{hyperref}

\newtheorem{theorem}{Theorem}
\newtheorem{proposition}{Proposition}
\newtheorem{conjecture}{Conjecture}
\newtheorem{corollary}{Corollary}

\begin{document}
\title{Using the Jones Polynomial to Prove Infinite Families of Knots Satisfy the Cosmetic Surgery Conjecture}
\author{F. M. Brady}
\date{December 23, 2025}
\maketitle
\abstract{This paper computes the Jones polynomial and the invariants obstructing cosmetic surgery which are derived from it for two infinite families of knots, proving they satisfy the Purely Cosmetic Surgery Conjecture. Both the method of computation and the method for generating families of knots extend.}

\section{Introduction}
Let $S^3_r(K)$ be the manifold obtained from surgery along $K$ with slope $r$.

\begin{conjecture}[Purely Cosmetic Surgery Conjecture] If $K$ is a knot and $S^3_{r}(K) \cong S^3_{r'}(K)$ via an orientation preserving homeomorphism then either $K$ is the unknot or $r = r'$.
\end{conjecture}

For a knot $K$ to have cosmetic surgeries, the surgery slopes $r$ and $r'$ must be $\pm 2$, the genus of $K$ must be $2$, and the Alexander polynomial of $K$ must be identically $1$ \cite{DEL}. The knot $K$ must also be prime \cite{Tao} and its Jones polynomial $V_K$ must satisfy $V_K''(1) = 0$ and $V_K'''(1) = 0$. \cite{IchiharaWu}. The derivatives of Jones polynomials are examples of finite type invariants; the higher order derivatives are powerful but previous obstructions were more complicated relations dependent on surgery slopes \cite{Ito}. We will derive, from \cite{Ito} and the restrictions on the Alexander polynomial from \cite{DEL}, the following.

\begin{proposition}
If a knot $K$ has cosmetic surgeries, then $V_K''''(1) = 0$.
\end{proposition}

This paper presents a method for computing the Jones polynomial, and, in particular, the derivatives of the Jones polynomial evaluated at $1$, for infinite families of knots obtained from additions of twists at crossings. The method enables direct computation of arbitrarily high derivatives of the Jones polynomial for these infinite families, allowing higher order derivative restrictions to be computationally useful in general.

\begin{theorem} The infinite families of knots generated from the knots $7_6$, $8_{12}$ and $10_{58}$ satisfy the Cosmetic Surgery Conjecture.
\end{theorem}

The paper begins by analyzing previous work on the conjecture, then proving the above proposition, and establishing the computational method for the Jones polynomial. It then considers specifically the case of $7_6$, first deriving the base case information needed to apply the method, then working out the casework to check the entire family and verifying that the cases which are not able to be disproven are in fact unknots in disguise. Finally, the case of $10_{58}$ is treated in the same way, with the omission of the unknot checking, as there are no unknots in this family. The case of $8_{12}$ is proven as a corollary of the Alexander polynomial computations for $10_{58}$.\\

\subsection{Acknowledgements} I would like to thank Zoltán Szabó for suggesting I investigate this question and for many helpful conversations. I would also like to thank Peter Ozsváth for his support.

\section{Relation to Previous Work}
\begin{figure}
\begin{center}
\includegraphics[width=.5\linewidth]{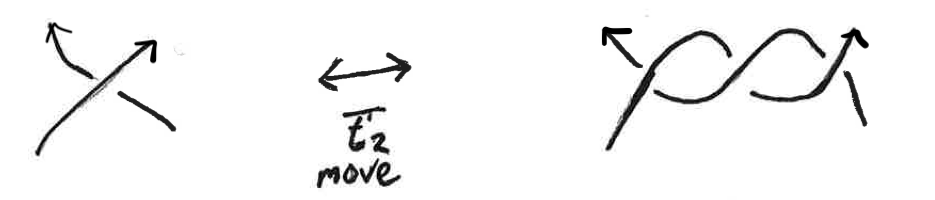}
\end{center}
\caption{A $\overline{t_2'}$ move. Every $\overline{t_2'}$ move is of this form or its mirror is.}
\label{t2-move}
\end{figure}

The knots considered here are knots of canonical genus $2$. The canonical genus of a knot is the minimal genus obtained from the Seifert algorithm over any possible diagram. While the potential difference between the genus and the canonical genus of any knot is unbounded \cite{KobayashiKobayashi}, considering knots of a given canonical genus is a useful way of ordering families of knots. Explicitly, all knots of canonical genus $2$ can be related by $\overline{t_2'}$ moves to one of the knots: $5_1$, $6_2$, $6_3$, $7_5$, $7_6$, $7_7$, $8_{12}$, $8_{14}$, $8_{15}$, $9_{23}$, $9_{25}$, $9_{38}$, $9_{39}$, $9_{41}$, $10_{58}$, $10_{97}$, $10_{101}$, $10_{120}$, $11_{123}$, $11_{148}$, $11_{329}$, $12_{1097}$, $12_{1202}$, or $13_{4233}$ \cite{Stoimenow}. A $\overline{t_2'}$ move is the operation of changing a single crossing to three crossings, which corresponds to introducing a full twist to a band in a Seifert diagram, visualized as a collection of nodes and connecting bands. See \ref{t2-move}.

There is some overlap between canonical genus $2$ families of knots and knots previously shown to satisfy the Cosmetic Surgery Conjecture, but enumerating these overlaps verifies that the result is novel. For instance, the knots related to $5_1$ are precisely the genus $2$, $5$-stranded pretzel knots. Pretzel knots have already been proven to satisfy the conjecture \cite{StipsiczSzabo}.

All prime alternating knots of genus $2$ and signature $0$ (hence $\tau$ invariant $0$) can be obtained by starting with one of six knots, $6_3$, $7_7$, $8_{12}$, $9_{41}$, $10_{58}$, $12_{1202}$, and performing $\overline{t_2'}$ moves. (Proposition 3.2; \cite{Stoimenow}) These are not excluded from consideration entirely, since it is possible to introduce crossing changes that do not preserve the alternating structure, but for part of the family the result has already been shown; \cite{IchiharaJong} shows the result for a subset of these and \cite{DEL} proves the result for alternating knots in general.

Previous work has proven many two bridge knots satisfy the conjecture \cite{ichiharaSaito}. As all two bridge knots are alternating \cite{Goodrick}, there may be some overlap with each of the families of canonical genus $2$ knots, but these are limited as the particular collection of signs must begin as an alternating knot, and most of the choices of sign do not give alternating knots. Additionally, since the bridge number cannot decrease, when starting with a knot for which the bridge number is greater than $2$, the entire family will be disjoint.

The overlap with $3$ braid knots, which have been proven to satisfy the conjecture \cite{Varvarezos}, is finite, since if a $3$ braid knot has genus $2$ it must have $10$ or fewer crossings, so only a finite number of cases from each knot could have possible overlap.

Previous work \cite{IchiharaJong} has shown some Montesinos knots satisfy the conjecture. According to the classifications of the paper, using their notation for Montesinos knots, all of their classes but (e2)-(e3) can be seen to be canonical genus $2$ knots, as shown below. However, there are many infinite families of canonical genus $2$ knots that are not Montesinos knots so this category is more general.

\begin{proposition} The Montesinos knots not proven to lack cosmetic surgeries by \cite{IchiharaJong} are a subset of the families generated by $6_2$, $7_5$, $7_6$, $7_7$ and $10_{58}$.
\end{proposition}

Here, $n_x$ is the signed number so that $(2 |n_x| - 1) = |2x + 1|$ and $n_x$ has the same sign as $2x + 1$, and all variables are integers with restrictions as indicated. Going forwards, the diagrams for knot skeletons will be nodes, represented by loops and bands, represented either by a single line, in the case of odd twists, or a double line, in the case of even twists.

\begin{figure}
\begin{center}
\includegraphics[width=.4\linewidth]{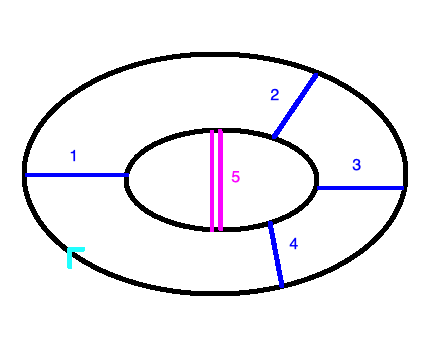}
\includegraphics[width=.4\linewidth]{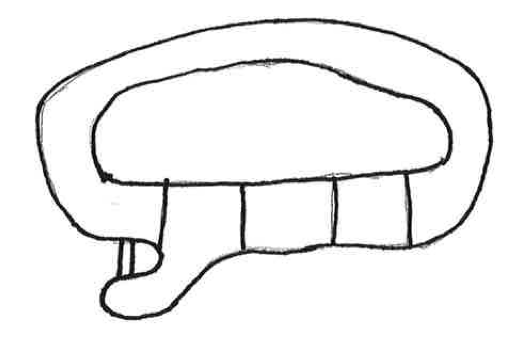}
\end{center}
\caption{The skeletons for the knot family built from $6_2$ and (o1) Montesinos knots.}
\label{o1}
\end{figure}

(o1), $a, c, d \neq -1$, and $a,b,c,d \neq 0$. Proven by \cite{IchiharaJong} to satisfy the conjecture for (o1') where additionally $e$ is either $0$ or $-1$.
$$M([2a+1,2b],[2c+1],[2d+1],[2e+1])= 6_2(n_a, n_c, n_d, n_e, b)$$ 
See figure \ref{o1} for comparison.

\begin{figure}
\begin{center}
\includegraphics[width=.3\linewidth]{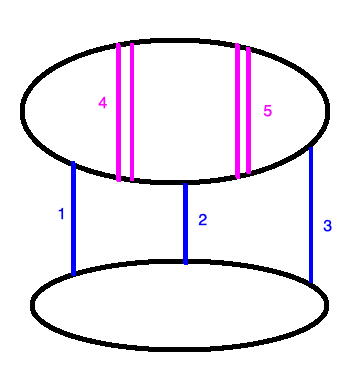}
\includegraphics[width=.4\linewidth]{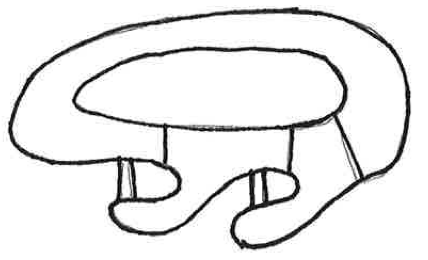}
\end{center}
\caption{The skeletons for the knot family built from $7_7$ and (o2) Montesinos knots.}
\label{o2}
\end{figure}

(o2): $a, c, d \neq 0, -1$, $b, d \neq 0$.
\[ M([2a+1, 2b], [2c+1, 2d], [2e+1]) = 7_7(n_a, n_e, n_c, b, d) \]
See figure \ref{o2} for comparison.

\begin{figure}
\begin{center}
\includegraphics[width=.2\linewidth]{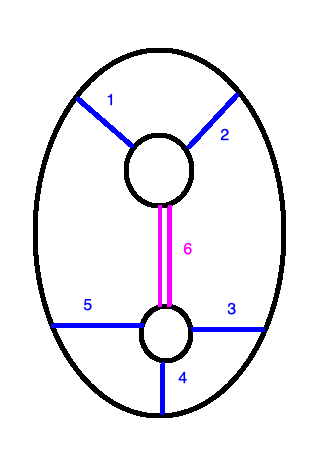}
\includegraphics[width=.3\linewidth]{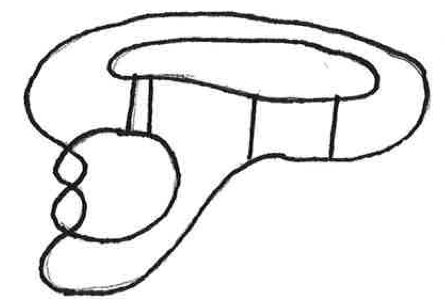}
\includegraphics[width=.3\linewidth]{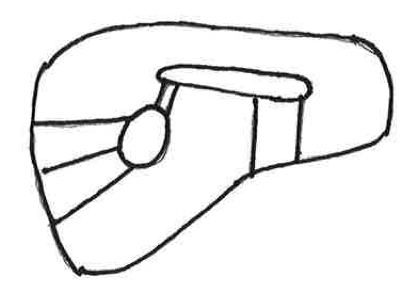}
\end{center}
\caption{The skeletons for the knot family built from $7_5$ and (o3) Montesinos knots, and the generalization of (o3) knots which connects the two.}
\label{o3}
\end{figure}

(o3), $b, c \neq -1, 0$, $a \neq 0$. Proven by \cite{IchiharaJong} to satisfy the conjecture for (o3'), where additionally $|a| = 1$. 
\[ M([2a, \pm 3], [2b + 1], [2c + 1]) = 7_5(n_c, n_b, \pm 1, \pm 1, \pm 1, a) \]
See figure \ref{o3} for comparison.

\begin{figure}
\begin{center}
\includegraphics[width=.2\linewidth]{7_5-skeleton.png}
\includegraphics[width=.3\linewidth]{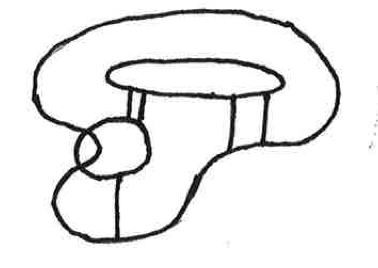}
\includegraphics[width=.3\linewidth]{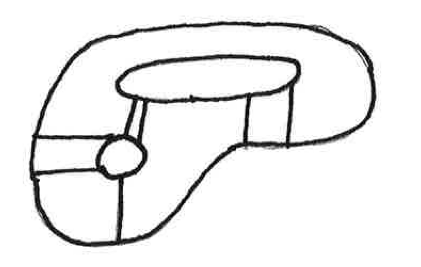}
\end{center}
\caption{The skeletons for the knot family built from $7_5$ and (o4) Montesinos knots, and the generalization of (o4) knots which connects the two.}
\label{o4}
\end{figure}

(o4) $b, c, d \neq -1, 0$, $a \neq 0$. Proven by \cite{IchiharaJong} to satisfy the conjecture for (o4'), where $|a| = 1$. 
\[ M([2a, \pm 2, 2b + 1], [2c + 1], [2d + 1]) = 7_5(n_d, n_c, n_b, \pm 1, \pm 1, a)\]
See figure \ref{o4} for comparison.

\begin{figure}
\begin{center}
\includegraphics[width=.3\linewidth]{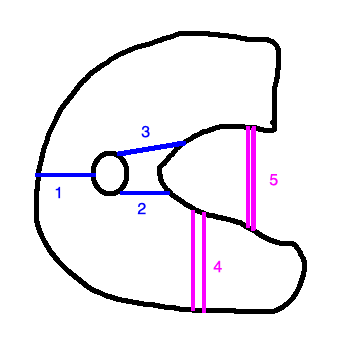}
\includegraphics[width=.4\linewidth]{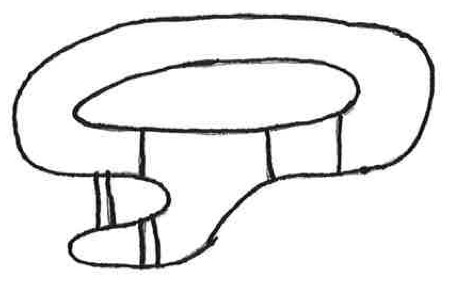}
\end{center}
\caption{The skeletons for the knot family built from $7_6$ and (o5) Montesinos knots.}
\label{o5}
\end{figure}

(o5): $a, d, e \neq -1, 0$, $b, c \neq 0$.
$$M([2a+1, 2b, 2c], [2d+1], [2e+1]) = 7_6(n_a, n_d, n_e, c, b)$$
See figure \ref{o5} for comparison.

\begin{figure}
\begin{center}
\includegraphics[width=.3\linewidth]{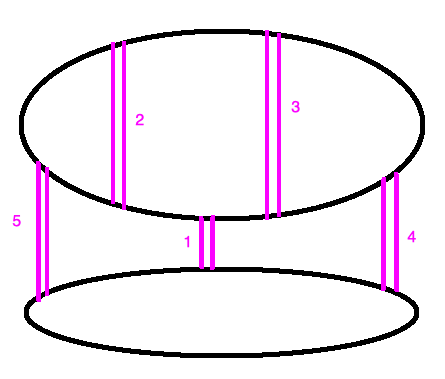}
\includegraphics[width=.4\linewidth]{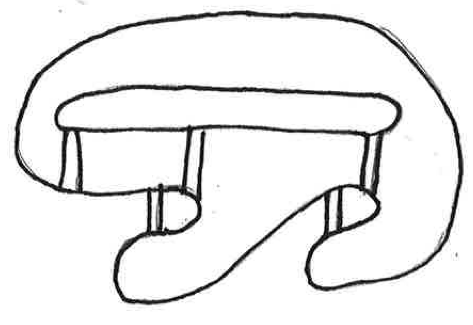}
\end{center}
\caption{The skeletons for the knot family built from $7_6$ and (e1) Montesinos knots.}
\label{e1}
\end{figure}

(e1): $a, b, c, d, e \neq 0$
$$M([2a], [2b, 2c], [2d, 2e]) = 10_{58}(a, c, e, d, b)$$
See figure \ref{e1} for comparison.

\begin{figure}
\begin{center}
\includegraphics[width=.4\linewidth]{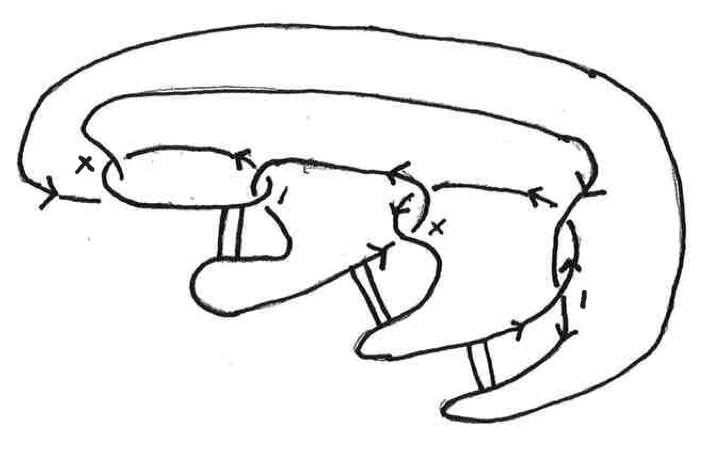}
\end{center}
\caption{The skeleton for (e2) Montesinos knots.}
\label{e2}
\end{figure}

(e2): $a, b, c \neq 0$. $$M([2],[-2, 2a], [2, 2b], [-2, 2c])$$ is not of canonical genus $2$ and has been proven to satisfy the conjecture \cite{IchiharaJong}. See figure \ref{e2}. Comparison with every knot family of canonical genus $2$ does not give any matches. Regardless, if there was overlap, counting parameters shows that (e2) would be a proper subset of any knot family of canonical genus $2$, as the smallest number of parameters that any family of canonical genus $2$ has is $4$.

\begin{figure}
\begin{center}
\includegraphics[width=.4\linewidth]{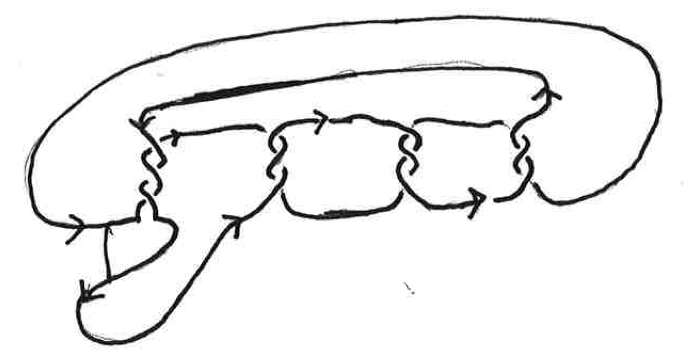}
\end{center}
\caption{The skeleton for (e3) Montesinos knots.}
\label{e3}
\end{figure}

(e3): $a \neq 0, -1$. $$M([3, 2a+1], [-3], [3], [-3])$$ is not of canonical genus $2$ and has been proven to satisfy the conjecture  \cite{IchiharaJong}. See figure \ref{e3}. Comparison gives no matches, but again, if there was overlap, it would be at most a proper subset of any canonical genus $2$ family, as it only has one parameter.

The families (o1), (o2), (o3), (o4), (o5) and (e1) are the remaining Montesinos knots which could possibly have cosmetic surgeries, and in fact, (o5) is almost exactly the family generated by $7_6$, with some differences on starting indices, and (e1) is precisely the family generated by $10_{58}$.

%\begin{figure}
%\includegraphics[width=.2\linewidth]{6_2-skeleton.png}
%\includegraphics[width=.2\linewidth]{7_7-skeleton.png}
%\includegraphics[width=.15\linewidth]{7_5-skeleton.png}
%\includegraphics[width=.2\linewidth]{7_6-skeleton.png}
%\includegraphics[width=.2\linewidth]{10_58-skeleton.png}
 %\caption{The skeletons for $6_2$, $7_7$, $7_5$, $7_6$, $10_{58}$, from left to right, with labeling to match}
%\end{figure}

\begin{corollary} If a Montesinos knot has a cosmetic surgery, it must be in one of the families (o1), (o2), (o3) or (o4).
\end{corollary} 

A final restriction can be obtained by evaluating the Jones polynomial itself at a fifth root of unity, $\xi_5$. Using the results of \cite{Detcherry}, if a knot $K$ has a purely cosmetic surgery the slope of the surgery will be of the form $\frac{1}{5k}$ unless $J_k(\xi_5) = 1$. Hence as $2$ is not of this form, if $J_K(\xi_5) \neq 1$, we may exclude the possibility of cosmetic surgery. In this context, this will give a mod $5$ congruence restriction on the number of pairs of twists in each possible location. Although this will not be needed in this computation, it works very well with these methods and is a potentially useful way of identifying infinitely many knots within a family which satisfy the conjecture.

These results build upon other foundational work, including \cite{GL} \cite{OS} \cite{Wang} \cite{Wu} \cite{NiWu} \cite{Hanselman} \cite{FPS}. Additionally, Ren has proven \cite{Ren} that if the Cosmetic Surgery Conjecture is false, there must exist a hyperbolic knot which has cosmetic surgeries.

\section{Jones Polynomial Restrictions}
Now, considering the result of \cite{Ito} we can consider the finite type invariant of one degree higher, which gives a relation between the fourth derivative of the Jones polynomial evaluated at $1$ and the Alexander polynomial, and use the results of \cite{DEL} to prove the position stated in the introductory section, that if $K$ has cosmetic surgeries, $V_K''''(1) = 0$.

\begin{proof} From Corollary 1.5(i) of \cite{Ito}, we have
\[
p^2(24w_4(K) - 5v_4(K)) + 5v_4(K) + q^2(210v_6(K) + 5v_4(K)) = 0
\]

We have that the slope is $\pm 2$, so $p^2 = 4$ and $q^2 = 1$, implying
\[
4(24w_4(K) - 5v_4(K)) + 5v_4(K) + (210v_6(K) + 5v_4(K)) = 0
\]
\[
96 w_4(K) + 210v_6(K) - 10v_4(K) = 0.
\]
Additionally, by Lemma 2.1 of \cite{Ito} we have
\[
v_4(K) = \frac{-1}{2} a_4(K) - \frac{1}{24} a_2(K) + \frac{1}{4} a_2(K)^2
\]
\[
w_4(K) = \frac{1}{96} j_4(K) + \frac{3}{32} a_4(K) - \frac{9}{2} a_2(K)^2
\]
\[
v_6(K) = \frac{-1}{2} a_6(K) - \frac{1}{12} a_4(K) - \frac{1}{720} a_2(K) + \frac{1}{24} a_2(K)^2 + \frac{1}{2} a_2(K) a_4(K) - \frac{1}{6} a_2(K)^3
\]
where $a_{2i}(K)$ is the coefficient of $z^{2 i}$ in the Conway polynomial and $j_n$ is the coefficient of $h^n$ in the Jones polynomial $V_K(e^h)$. However, if $K$ has a cosmetic surgery, the Alexander polynomial and hence the Conway polynomial must be trivial \cite{DEL}. Therefore, we can simplify these formulae to:
\[
v_4(K) = 0
\]
\[
w_4(K) = \frac{1}{96} j_4(K)
\]
\[
v_6(K) = 0.
\]

Plugging these back into the equation above, we have that $j_4(K) = 0$, so $V_K''''(1) = 0$ as well for a knot that satisfies the Cosmetic Surgery Conjecture.
\end{proof}

\section{General Jones Polynomial Computation}

\begin{theorem}
If $K(n_1, ... n_k)$ is the result of adding $s_i(2|n_i| - m_i)$ twists to the $i$th crossing of a knot $K$, then the Jones polynomial of this new knot is:
\[
V_{K(n_1, ... n_k)}(t) = \sum \limits_{x \in \{a,b\}^k} \left( \prod \limits_{i= 1}^k c^{s_i}(t, |n_i|, x) \right) V_x(t)
\]
where $c^{\pm}(t,n,0) = \frac{1-t^{\pm 2n}}{1-t^{\pm 2}} \left( \pm t^{\pm 1} \right) \left( t^{\frac{1}{2}} - t^{\frac{-1}{2}} \right) $ and $c^{\pm}(t,n,1) = t^{\pm 2n}$, and $V_x(t)$ is the Jones polynomial for the knot obtained by resolving every crossing of $K$ corresponding to a $0$ in the corresponding coordinate of $x$, and ensuring there are $-s_i m_i$ twists at the $i$th crossing for every $i$ where $x_i = 1$.
\end{theorem}

Note that in the above formula, there are a finite number of links whose Jones polynomials must be computed for any given $K$.

\begin{proof}
The first step is to find the Seifert surfaces as a set of twisted bands. The $\overline{t_2’}$ moves will correspond to adding pairs of additional twists in each band. These bands will provide the necessary structure to orderly analyze the construction.
To compute the Jones polynomial in general, consider a band with $\pm (2n-1)$ twists in the case of an odd twist or $\pm 2n$ in the case of an even twist. Consider any crossing in this band. See \ref{fig2} 
 \begin{figure}
\; \; \; \; \; \; \; \; \; \includegraphics[width=.5\linewidth]{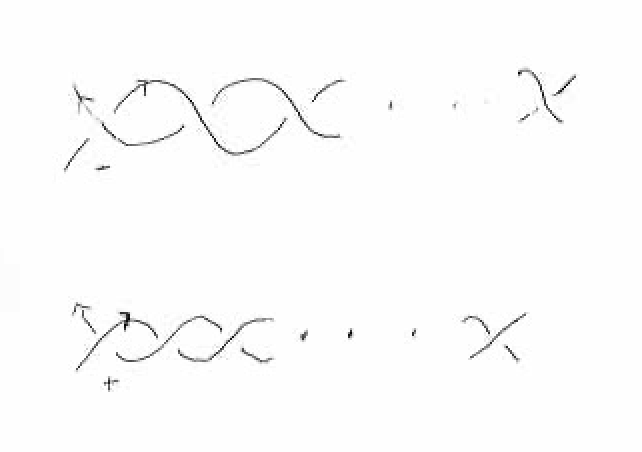}
 \caption{Crossings}
 \label{fig2}
\end{figure}

If the crossing is reversed, this will introduce a twist of opposite parity which, using the Reidemeister II move, can cancel one of the starting twists. Therefore, reversing the crossing sends the knot from one with $\pm(2n - m)$ twists on this band to $\pm(2(n - 1) - m)$ twists $(m \in \{0, 1 \})$.\

If the crossing is resolved, this will cut the band in half. The result of the untwisting, using repeated Reidemeister I moves, does not depend on the original $n$ (or $m$).

The skein relation for the Jones polynomial can thus connect these three quantities and create a recursion relation. This will end with two polynomials: the Jones polynomial for the link with $m$ twists on the band, and the Jones polynomial for the link with no band. Explicitly, let $Q(t)$ be the Jones polynomial for the knot (or link) with no band, and $P_n(t)$ be the polynomial for the knot (or link) with $\pm(2n - m)$ twists, for positive or negative twists. Further, differentiate these cases based on if the original twist was positive or negative.

For positive starting twists:
$$
\begin{array}{rcl} (t^{\frac{1}{2}} - t^{\frac{-1}{2}})Q^+(t) & = & t^{-1}P_n^+(t) - t P_{n - 1}^+(t) \\ t (t^{\frac{1}{2}} - t^{\frac{-1}{2}}) Q^+(t) + t^2 P_{n - 1}^+(t) & = & P_n^+(t) \\ (t + t^3) (t^{\frac{1}{2}} - t^{\frac{-1}{2}}) Q^+(t) + t^4 P_{n - 2}^+(t) & = & P_n^+(t) \\ (1 + t^2 + ... + \;t^{2(n - 1)}) t (t^{\frac{1}{2}} - t^{\frac{-1}{2}}) Q^+(t) + t^{2n} P_{0}^+(t) & = & P_n^+(t) \\ \frac{1 - t^{2n}}{1 - t^2} t (t^{\frac{1}{2}} - t^{\frac{-1}{2}}) Q^+(t) + t^{2n} P_{0}^+(t) & = & P_n^+(t) \\ \end{array}
$$

For negative starting twists:
$$
\begin{array}{rcl} (t^{\frac{1}{2}} - t^{\frac{-1}{2}})Q^-(t) & = & t^{-1}P_{n - 1}^-(t) - t P_{n}^-(t) \\ (-1)(t^{-1}) (t^{\frac{1}{2}} - t^{\frac{-1}{2}}) Q^-(t) + t^{-2} P_{n - 1}^-(t) & = & P_n^-(t) \\ (1 + t^{-2}) (-t^{-1})(t^{\frac{1}{2}} - t^{\frac{-1}{2}}) Q^-(t) + t^{-4} P_{n - 2}^-(t) & = & P_n^-(t) \\ \frac{1 - t^{-2n}}{1 - t^{-2}} (-t^{-1}) (t^{\frac{1}{2}} - t^{\frac{-1}{2}}) Q^-(t) + t^{-2n} P_{0}^-(t) & = & P_n^-(t) \\ \end{array}
$$

Therefore,
$$
P_n^{\pm}(t) = t^{\pm 2n} P_0^{\pm}(t) + \frac{1 - t^{\pm 2n}}{1 - t^{\pm 2}}(\pm t^{\pm 1}) (t^{\frac{1}{2}} - t^{\frac{-1}{2}}) Q^{\pm}(t)
$$
where $P_0^{\pm}(t)$ and $Q^{\pm}(t)$ will depend on the specific knot. 

Using the fact that the skein relation can be applied to different parts of a diagram in succession, and it is by nature local (while further simplification may be possible in a given case, it is by no means required), this process can be repeated once for each band, with each final link term picking up one of the prefactors $t^{\pm 2n_i}$ or $\frac{1 - t^{\pm 2n_i}}{1 - t^{\pm 2}} (\pm t^{\pm 1}) (t^{\frac{1}{2}} - t^{\frac{-1}{2}})$ depending on whether a band remains or not in the $i$th case, where there began with $\pm (2n_i - m_i)$ twists, and ending with either no band or a band with $\mp m_i$ twists (so either $0$ in the case of an even band or, $1$ or $-1$ in the case of an odd band). Repeating this analysis for every band will give the desired formula. 
\end{proof}

To compute the derivatives of the Jones polynomial at $1$, it is necessary to evaluate the $k$th derivatives of each of the prefactors, evaluated at $1$. For $k = 0, .. 4$, these are as follows.
$$
\begin{array}{c||c|c|c|c} k & c^{\pm}(t,n,0) & c^{\pm}(t,n,1) \\ \hline 0 & 0 & 1\\ 1 & \pm n & \pm 2n\\ 2 & 2n^2 \mp n & 2n (2n \mp 1) \\ 3 & \frac{n}{4}(\pm 5 - 24n \pm 16n^2) & 4n \left( \pm 1 - 3n \pm 2 n^2 \right) \\ 4 & \frac{n}{2} (\mp 3 + 38 n \mp 48 n^2 + 16 n^3) & 4 n (\mp 3 + 11 n \mp 12 n^2 + 4 n^3) \\ \end{array}
$$
For the explicit computation, is simple to use these values to compute the derivative of a sum of products numerically.

\section{Base Cases for $7_6$}
\begin{figure}
\begin{center}
\includegraphics[width=.3\linewidth]{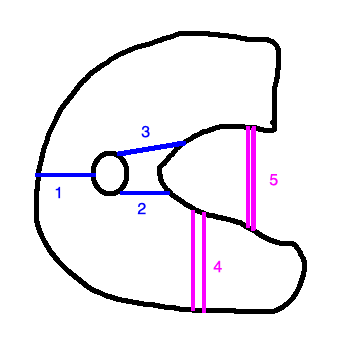}
\end{center}
 \caption{$7_6$ Skeleton}
 \label{fig1}
\end{figure}

We can construct the Seifert surface according to the Seifert algorithm for $7_6$ as in \ref{fig1}. The labeling on the bands will be arranged according to the labeled numbers. As before, a double line will denote an evenly twisted band, and a single line will denote an oddly twisted one.

Despite a seeming preponderance of base cases, patterns emerge which reduces the number of explicit calculations needed. In the case of $7_6$ analyzing the even bands allows us to find patterns to reduce the computation required. Using also that the Jones polynomial is multiplicative under connected sum (and does not detect the different knots formed by the choices available when performing a connected sum of links)
$$
\begin{array}{rcl}
7_6(x_1, x_2, x_3, 0, 0) & = & X_3(x_1,x_2,x_3) \\
7_6(x_1, x_2, x_3, 1, 0) & = & X_3(x_1,x_2,x_3) \sqcup U  \\
7_6(x_1, x_2, x_3, 0, 1) & = & X_2(x_1,x_2) \# X_1(x_3) \\
7_6(x_1, x_2, x_3, 1, 1) & = & X_3(x_1,x_2,x_3) \\
\end{array}
$$
Disjoint union with the unknot will result in a multiplication by $-\left(t^{\frac{1}{2}} + t^{\frac{-1}{2}} \right)$. $X_n$ is the pattern of $n$ odd twist bands connecting two disjoint circles.

To compute $X_n$ note that any inputs to $X_n$ which are $0$s can be removed and the $n$ reduced correspondingly. Additionally, $+$ and $-$ twists adjacent to each other cancel out and can be removed, with the $n$ reduced by $2$. Therefore, the only computations needed are $X_n(+, ..., +)$, $X_n(-, ..., -)$ for general $n$. Additionally, using the skein relation at a given crossing gives a relation between $X_n$ and $X_{n - 1}$, where the crossing is removed and $X_{n - 2}$, when the crossing is flipped, so
$$
\begin{array}{r|l}
X_0() & 2U \\
X_1(\pm) & U \\
X_n(+, ... +) & t^2 X_{n - 2}(+, ..., +) + \left(t^{\frac{3}{2}} - t^{\frac{1}{2}} \right) X_{n - 1}(+, ..., +)\\
X_n(-, ... -) & t^{-2} X_{n - 2}(-, ..., -) + \left( t^{\frac{-3}{2}} - t^{\frac{-1}{2}}\right) X_{n - 1}(-, ..., -) \\
\end{array}
$$
Combining all of these facts gives a concise way of expressing all of the base cases for computation of the Jones polynomial of every knot formed from $7_6$.

\section{Alexander Polynomial Computation for $7_6$}
Similarly, we can compute the Alexander polynomial by constructing the Seifert matrix. The main result is that for a band with $n$ twists (positive or negative), the pushoff of a curve running along it will have $-n$ signed self-intersections along the curve. Counting the self-intersections will give twice the linking number. For even numbers of twists, these will be all of the self-intersections. For odd numbers of twists, there will be an additional crossing so as to get the pushed off band on the same side as it began (relative to the given projection).

For $7_6$ we will find the matrix relative to the below curves and standard right hand orientation of the Seifert surface given the orientation of the boundary as shown. The representative curves will run along the $1$st and $2$nd bands ($x_1$), the $1$st and $3$rd bands ($x_2$), the $4$th band ($x_3$), and the $5$th band ($x_4$).

Using that there are $n_i$ twists along the $i$th band (using labeling as above), all curves give self-linking numbers of:
$$
S_{7_6}(n_1,n_2,n_3,n_4,n_5) = \left( \begin{array}{cccc} \frac{-(n_1 + n_2)}{2} & & & \\ & \frac{-(n_1 + n_3)}{2} & & \\ & & \frac{-n_4}{2} & \\ & & & \frac{-n_5}{2} \\ \end{array} \right)
$$
The pushoff of $x_4$ is disjoint from all other curves, but the pushoffs of the other curves go below $x_4$ giving two positive crossings. The curve $x_3$ does not interact with either $x_1$ or $x_2$.  
$$
S_{7_6}(n_1,n_2,n_3,n_4,n_5) =\left( \begin{array}{cccc} \frac{-(n_1 + n_2)}{2} & & 0 & 1\\ & \frac{-(n_1 + n_3)}{2} & 0 & 1\\ 0 & 0 & \frac{-n_4}{2} & 1\\ 0& 0& 0& \frac{-n_5}{2} \\ \end{array} \right)
$$
Both the pushoff of $x_1$ and the pushoff of $x_2$ have $n_1$ crossings with the other along the first band. The pushoff of each goes above the other along the innermost Seifert surface disk. 

This gives a positive crossing in the pushoff of $x_2$ and a negative crossing in the pushoff of $x_1$, giving final Seifert matrix of:
$$
S_{7_6}(n_1,n_2,n_3,n_4,n_5) = \left( \begin{array}{cccc} \frac{-(n_1 + n_2)}{2} & \frac{-n_1 - 1}{2} & 0 & 1\\ \frac{-n_1+1}{2} & \frac{-(n_1 + n_3)}{2} & 0 & 1\\ 0 & 0 & \frac{-n_4}{2} & 1\\ 0& 0& 0& \frac{-n_5}{2} \\ \end{array} \right)
$$
To translate this into the language that will be used to compute things, if the sign of the $i$th crossing is $s_i$ and the number of total twists applied are $t_i$ we have for an even band $n_i = 2 s_i t_i$ and for an odd band $n_i = 2s_i n_i - s_i$. Therefore, this becomes:
$$
\left( \begin{array}{cccc} -s_1 t_1 - s_2 t_2 + \frac{s_1 + s_2}{2} & -s_1 t_1 + \frac{s_1 - 1}{2} & 0 & 1\\ -s_1 t_1 + \frac{s_1+1}{2} & -s_1 t_1 - s_2 t_2 + \frac{s_1 + s_3}{2} & 0 & 1\\ 0 & 0 & -s_4 t_4 & 1\\ 0& 0& 0& -s_5 t_5 \\ \end{array} \right)
$$

\section{Computations for the Knot Family Generated by $7_6$}
For the knot $7_6$, there are $2^5 = 32$ cases for signs. All will be handled separately, both because the formulas change and because this ensures that $t_i$ is a positive integer in every case (which will enable showing specific quantities must be strictly positive or negative and cannot be zero).

The cases are as follows:\\

$(+++++)$, $(+++--)$, $(---++)$, $(-----)$\\
The leading term of the Alexander polynomial is 
\[ ((b - 1) (c - 1) + a (b + c -1)) d e,\]
which is strictly positive as it is the sum of a nonnegative and a positive number, hence non-vanishing.\\

$(++++-)$, $(+++-+)$, $(---+-)$, $(----+)$\\
The leading term of the Alexander polynomial is 
\[ -((b - 1) (c - 1) + a (b + c -1)) d e,\] 
which is strictly negative as it is the sum of a nonpositive and a negative number, hence non-vanishing.\\

$(++-++)$,$(--+--)$\\
The leading term of the Alexander polynomial is $(ab - c(a+b-1))de$. If this is zero, we have that $c = \frac{a b}{a + b - 1}$. The second derivative of the Jones polynomial evaluated at $1$ is:
$$
V_K''(1) = -6(ab - c(a+b+d-1)+ d(b+e))
$$
Substituting for $c$ gives $V_K''(1)$ becomes:
$$
V_K''(1) = \frac{-6 d (b^2 + b (e-1) + (a-1) e)}{a + b -1}
$$
This is negative, so this cannot vanish and we can exclude this case.\\

$(++--+)$, $(--++-)$\\
The leading term of the Alexander polynomial is $(-ab + c(a+b-1))de$. If this is zero, we have $c = \frac{a b}{a + b - 1}$. $V_K''(1)$ is $-6(ab + c(1-a-b+d)-d(b+e))$. Substituting for $c$ gives $V_K''(1)$ becomes $\frac{6 d (b^2 + b (e-1) + (a-1) e)}{a + b - 1}$, which is always positive.\\ 

$(+-+++)$, $(-+---)$\\
The leading term of the Alexander polynomial is $(ac - b(a + c - 1))de$. If this vanishes $b = \frac{a c}{a + c - 1}$. $V_K''(1)$ is $-6(ac + d(c+e)-b(a+c+d-1))$ so substituting for $b$ gives that this becomes $\frac{-6 d (c^2 + c (e-1) + (a-1) e)}{a + c-1}$, which is always negative, so it is impossible for both to vanish simultaneously.\\

$(+-+-+)$, $(-+-+-)$\\
The leading term of the Alexander polynomial is $(-ac + b(a+c-1))de$. If this vanishes $b = \frac{a c}{a + c - 1}$. $V_K''(1)$ is $6(-ac + b(a + c - d - 1) + d(c+e))$. If we substitute for $b$, we get that $V_K''(1)$ becomes $\frac{6 d (c^2 + c (e-1) + (a-1) e)}{a + c-1}$ which is always positive, so both quantities cannot vanish simultaneously.\\

$(+--+-)$, $(-++-+)$\\
The leading term of the Alexander polynomial is $(-bc + a(b+c-1))de$. If this vanishes, $a = \frac{b c}{b + c - 1}$. $V_K''(1)$ is $6(-bc + a(b+c-1) + d(b+c+e-1))$. Substituting for $a$ gives that $V_K''(1)$ becomes: $6d(b+c+e-1)$. This is immediately always positive, so both cannot vanish simultaneously.\\

$(+----)$, $(-++++)$\\
The leading term of the Alexander polynomial is $(bc - a(b + c - 1)) de$. If this vanishes, $a = \frac{b c}{b+c-1}$. $V_K''(1)$ is $6(-bc + a(b+c-1) - d(b + c + e - 1))$. Substituting for $a$ gives that $V_K''(1)$ becomes $-6d(b+c+e-1)$. This is always negative, so this cannot vanish when the first term of the Alexander polynomial does.\\

Now we get to the cases where some of the knots in question are unknots. \\

$(++-+-)$\\
The leading term of the Alexander polynomial is $(-ab + c(a+b-1)) de$. If this is zero, we have $c = \frac{a b}{a + b - 1}$. $V_K''(1)$ is $-6(ab - c(a+b-1)+d(b-e-c))$. Substituting for $c$ gives that $V_K''(1)$ becomes $\frac{6 d (-b^2 + (a-1) e + b (1 + e))}{a + b-1}$. Therefore this vanishes only if $-b^2 + (a-1) e + b (1 + e)$ does, which occurs only when $e = \frac{b (b -1)}{a + b - 1}$. Note that this implies $b \neq 1$ as $e$ cannot be zero. Simplifying $V_K'''(1)$ with these relations:
\[
V_K'''(1) \mapsto \frac{18 (a-1) a (b-1) b}{a + b-1}
\]
Because $b > 1$, this vanishes only if $a = 1$. This implies that $c = \frac{b}{b} = 1$, and $e = \frac{b(b - 1)}{b} = b - 1$. Therefore, the knot is of the form $(1, e+1,-1,d,-e)$. Since this is the only case when $V_K'''(1)$ vanishes and this is a form of the unknot, as will be seen below, this case can also be excluded from possible counterexamples.\\
Related: $(--+-+)$\\
The only place where the calculation departs from the above case is in the simplification of $V_K'''(1)$. Here
$$
V_K'''(1) \mapsto \frac{-18 (a-1) a (b-1) b}{a + b-1}
$$
Therefore, the computation is precisely the same and only knots not excluded are of the form $(-1, -(e+1), 1, -d, e)$, which are precisely the mirrors of the above case.\\

$(++---)$\\
The leading term of the Alexander polynomial is $(ab - c(a+b-1)) d e$. If this is zero, we have $c = \frac{a b}{a+b-1}$. $V_K''(1)$ is $-6(ab - c(a+b-1)+d(-b+e+c))$. Substituting for $c$ gives that $V_K''(1)$ becomes $\frac{6 d (b^2 + e - a e - b (1 + e))}{a + b-1}$. Therefore this vanishes only if $e = \frac{b(b-1)}{a+b-1}$. Note that this implies $b \neq 1$ as $e$ cannot be zero. Simplifying $V_K'''(1)$ with these relations gives
\[
V_K'''(1) \mapsto \frac{18 (a - 1) a (b-1) b}{a + b-1}
\]
This vanishes only when $a = 1$. In that case, $c = 1$ and $e = b - 1$. Therefore, the only case not excluded by this calculation is $(1,e+1,-1,-d,-e)$.\\
Related: $(--+++)$\\
The only difference from the above calculation is that $V_K'''(1)$ simplifies as
$$
V_K'''(1) \mapsto \frac{-18 (a - 1) a (b-1) b}{a + b-1}
$$
Therefore, everything is the same so the knot must be of the form $(-1,-(e+1),1,d,e)$. This is the mirror of the knot family above, so if those are all the unknot, these must all be the unknot as well.\\

$(+-++-)$\\
The leading term of the Alexander polynomial is $(-ac + b(a + c - 1)) de$. If this is zero, then $b = \frac{ac}{a + c - 1}$. $V_K''(1)$ is $6(-ac + b(a + c + d - 1) + d(e - c))$. Substituting for $b$ gives $\frac{6 d (-c^2 + (a-1) e + c (1 + e))}{a + c-1}$. This vanishes only if $e = \frac{c(c-1)}{a + c - 1}$. This means that $c$ cannot be $1$ as $e$ cannot be $0$. Substituting for both $b$ and $e$ gives
$$
V_K'''(1) \mapsto \frac{18 (a-1) a (c-1) c}{a + c-1}
$$
This vanishes only if $a = 1$. In that case, $b=1$ and $e = -1 + c$. Therefore, the only knots not excluded are of the form $(1, -1,e+1,d,-e)$. This will be shown to be the unknot below.\\
Related: $(-+--+)$\\
The only difference in the above case is that $V_K'''(1)$ simplifies to:
$$
V_K'''(1) \mapsto \frac{-18 (a-1) a (c-1) c}{a + c-1}
$$
The computation is the same so the knots must be of the form $(-1, 1,-(e + 1), -d, e)$. These are the mirrors of the knots above, so if the knots from the prior case are unknots, so are these. \\

$(+-+--)$\\
The leading term of the Alexander polynomial is $(a c - b (a + c-1)) d e$. This vanishes only if $b = \frac{a c}{a + c - 1}$. $V_K''(1)$ is $6 (-a c + b (-1 + a + c - d) + d (c - e))$. Substituting for $b$ gives $V_K''(1)$ becomes $\frac{6 d (c^2 + e - a e - c (1 + e))}{a + c - 1}$ so this vanishes only if $e = \frac{c(c-1)}{a+c-1}$. Using these relations $V_K'''(1)$ becomes:
$$
V_K'''(1) \mapsto \frac{18 (a-1) a (c-1) c}{a + c-1}
$$
Because $e$ cannot be zero, $c$ cannot be $1$, so this vanishes only if $a = 1$. In that case $b = 1$ and $e = c - 1$. Therefore, the only knots not excluded are of the form $(1,-1,e+1,-d,-e)$. These will be shown to be the unknot below. \\
Related: $(-+-++)$\\
The only difference in the above case is that $V_K'''(1)$ simplifies to:
$$
V_K'''(1) \mapsto \frac{-18 (a-1) a (c-1) c}{a + c-1}
$$
The calculation is the same, so the only knot not excluded is $(-1,1,-(e+1),d,e)$, which is again the mirror of the previous class of knots.\\

$(+--++)$\\
The leading term of the Alexander polynomial is $(bc - a(b + c - 1))de$. If this vanishes, $a = \frac{b c}{b + c - 1}$. Also, $V_K''(1) = 6(-bc + a(b + c - 1) + d(b+c-e-1))$. Substituting for $a$ gives $V_K''(1)$ becomes $6d((b + c - 1) - e)$ so $e = b + c - 1$. Substituting this into the Jones polynomial derivatives gives:
$$
V_K'''(1) \mapsto \frac{-18 (-1 + b) b (-1 + c) c}{b + c-1}
$$
Therefore, either $b = 1$ or $c = 1$. In both cases, $a = 1$, but in the first case $e = c$ and in the second $e = b$. Therefore, the only knot cases not explicitly ruled out are $(1,-1,-e,d,e)$ or $(1, -e,-1,d,e)$.\\
Related: $(-++--)$\\
The only difference in the previous analysis is that the third derivative of the Jones polynomial at $1$ simplifies to
$$
V_K'''(1) \mapsto \frac{18 (-1 + b) b (-1 + c) c}{b + c-1}
$$
The remaining analysis is the same so the knot cases not explicitly ruled out here are $(-1, 1, e,-d,-e)$ and $(-1, e, 1, -d,-e)$. \\

$(+---+)$\\
The leading term of the Alexander polynomial is $(-bc + a(b + c - 1)) de$, so if this vanishes, $a = \frac{b c}{b + c - 1}$ and $V_K''(1) = 6(a(b + c - 1) - bc + d(1-b-c+e))$, which, substituting for $a$, becomes $-6 d (-1 + b + c - e)$. If both of these are zero, we have additionally that $e = b + c - 1$. Substituting for this relation,
$$
V_K'''(1) \mapsto \frac{-18 (b-1) b (c-1) c}{b + c-1}
$$
Therefore, either $b = 1$ or $c = 1$. In either case, $a = 1$. In the first, $e = c$ and in the second $e = b$. Therefore, the only knots not excluded are $(1, -1, -e, -d, e)$ and $(1,-e,-1,-d, e)$. \\
Related: $(-+++-)$\\
The only difference with the previous case is that $V_K'''(1)$ simplifies as:
$$
V_K'''(1) \mapsto \frac{18 (b-1) b (c-1) c}{b + c-1}
$$
The remaining analysis is the same so the only two knot types not excluded are $(-1, 1, e, d,-e)$ and $(-1,e,1,d,-e)$.\\

\section{Showing the Remaining $7_6$ Knots are Unknots}
To consolidate, the knot exceptions are:
$$
\begin{array}{c|c|c|c|c}
Case & Knot & & Case & Knot \\ \hline
(++-+-) & (1,e+1,-1,d,-e) & & (--+-+) & (-1,-(e+1),1,-d,e) \\
(++---) & (1,e+1,-1,-d,-e) & & (--+++) & (-1,-(e+1),1,d,e) \\
(+-++-) & (1,-1,e+1,d,-e) & & (-+--+) & (-1,1,-(e+1),-d,e) \\
(+-+--) & (1,-1,e+1,-d,-e) & & (-+-++) & (-1,1,-(e+1),d,e) \\
(+--++) & (1,-1,-e,d,e) & & (-++--) & (-1,1,e,-d,-e) \\
(+---+) & (1,-1,-e,-d,e) & & (-+++-) & (-1,1,e,d,-e) \\
(+--++) & (1,-e,-1,d,e) & & (-++--) & (-1,e,1,-d,-e) \\
(+---+) & (1,-e,-1,-d,e) & & (-+++-) & (-1,e,1,d,-e) \\
\end{array}
$$

Note that switching total signs will result in a mirror of a knot, so if we show the knot exceptions above whose first twist is $+1$ are the unknot, this will show the second column as well. Additionally, as they are arranged, it implies and we will show that the twisting on the $d$ band is irrelevant. There are two cases that will be investigated, the bands $1$ and $3$ having twists $1$ and $-1$ respectively, and the bands $1$ and $2$ having twists $1$ and $-1$ respectively. These will address all cases in the first column, and the second column follows by considering the mirror images of these cases. 

\begin{figure}
  \includegraphics[width=.5\linewidth]{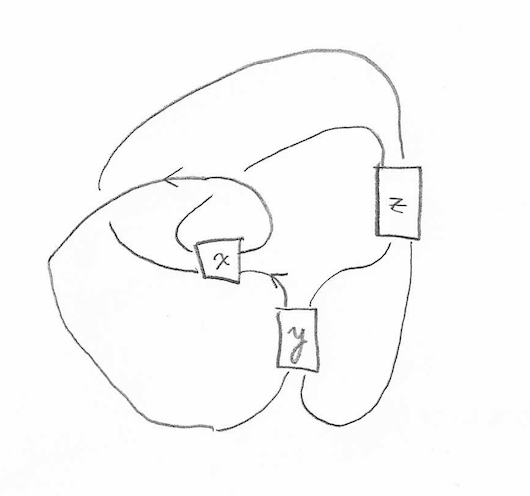}
  \includegraphics[width=.5\linewidth]{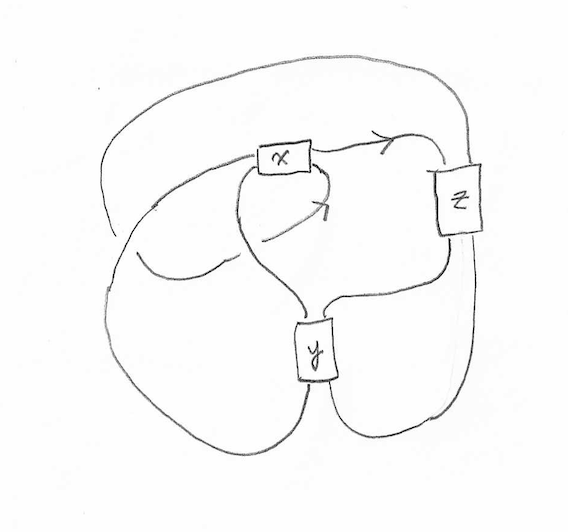}
  \caption{Cases A and B}
\end{figure}

\begin{figure}
\begin{center} \includegraphics[width=.3\linewidth]{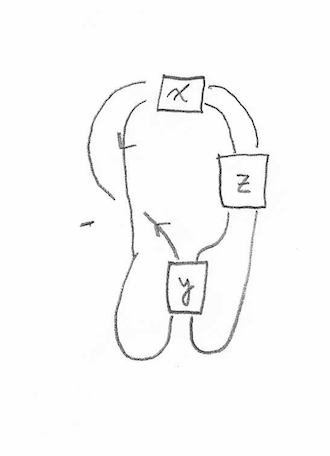} \; \; \; \; \; \;
  \includegraphics[width=.3\linewidth]{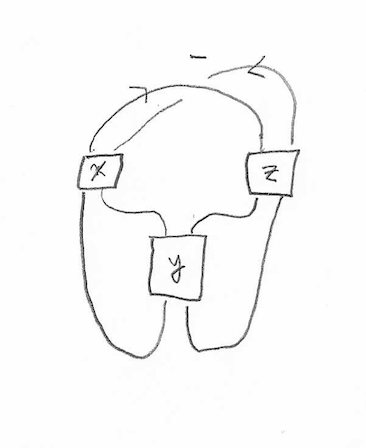}
  \end{center}
  \caption{Cases A and B simplified}
\end{figure}

\begin{figure}
\begin{center} \includegraphics[width=.3\linewidth]{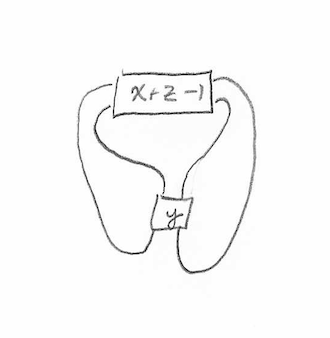} \; \; \; \; \; \; 
  \includegraphics[width=.2\linewidth]{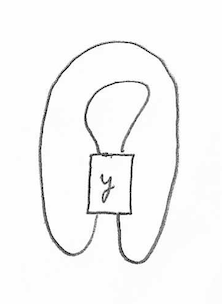}
  \end{center}
  \caption{Final Simplification of Both Cases}
\end{figure}

$$
\begin{array}{c|c|c}
Knot & Case & x + z - 1 \\ \hline
(1,e+1,-1,\pm d,-e) & A & (2(e+1)-1)+(-2 e) - 1 = 0 \\
(1,-e,-1, \pm d,e) & A & (2(-e)+1)+(2e) - 1 = 0 \\ 
 (1,-1,e+1, \pm d,-e) & B & (2(e+1)-1)+(-2e)-1 = 0 \\
 (1,-1,-e, \pm d,e) & B & (-2e+1)+(2e)-1 = 0 
\end{array}
$$

In all cases, the $y$ half twists can be readily untwisted to give the unknot, regardless of the value of $y$.

\section{Base Cases for $10_{58}$}
\begin{figure}
\begin{center}
\includegraphics[width=.4\linewidth]{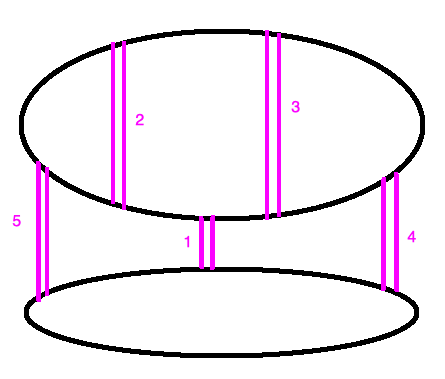}
\end{center}
\caption{$10_{58}$ skeleton.}
\label{10-58}
\end{figure}

Because all bands have even numbers of twists \ref{10-58}, all base cases will be disjoint unions of unknots. The Jones polynomial of a disjoint union of $n$ unknots is $\left(- \left(t^{\frac{1}{2}} + t^{\frac{-1}{2}} \right) \right)^{n - 1}$.
$$
\begin{array}{ccc|c}
n_5 & n_1 & n_4 & \text{ Number of Unknots } \\ \hline
0 & 0 & 0 & 2 + n_2 + n_3 \\
0 & 0 & 1 & 1 + n_2 + n_3 \\
0 & 1 & 0 & 1 + n_2 + n_3 \\
1 & 0 & 0 & 1 + n_2 + n_3 \\
0 & 1 & 1 & 2 + n_2 - n_3 \\
1 & 1 & 0 & 2 - n_2 + n_3 \\
1 & 0 & 1 & 2 - |n_3 - n_2| \\
1 & 1 & 1 & 3 - n_2 - n_3 \\
\end{array}
$$

\section{Alexander Polynomial Computations for $10_{58}$ and Analysis of $8_{12}$}
\begin{figure}
\begin{center}
\includegraphics[width=.4\linewidth]{10_58-skeleton.png}
\includegraphics[width=.4\linewidth]{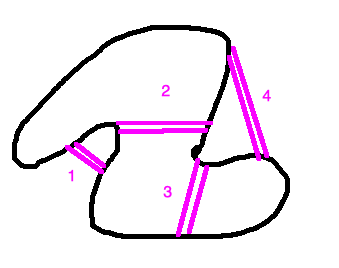}
\end{center}
\caption{Comparison of the $10_{58}$ and $8_{12}$ skeletons.}
\label{8-12}
\end{figure}

Consider bands passing through the $5$th and $1$st bands, the $1$st and $4$th bands, both going clockwise, one along the $2$nd band, and the $3$rd band. Then the self-linking numbers of the bands are:
$$
\left( \begin{array}{cccc} -n_1 s_1 - n_5 s_5 & & &\\
& -n_1 s_1 - n_4 s_4 & & \\
& & -n_2 s_2 & \\
&  & & -n_3 s_3 \\
\end{array} \right)
$$
The $3$rd and $4$th bands do not interact, and the $1$st and $2$nd have linking numbers of $n_1 s_1$ as their only interaction is going in different directions along the same band.
$$
\left( \begin{array}{cccc} -n_1 s_1 - n_5 s_5  & n_1 s_1 & &\\
n_1 s_1 & -n_1 s_1 - n_4 s_4 & & \\
& & -n_2 s_2 & 0 \\
& & 0 & -n_3 s_3 \\
\end{array} \right)
$$
The $1$st and $4$th bands do not interact, and the $2$nd and $3$rd do not either. The $1$st and $3$rd link with linking number $1$ only in one direction of pushoff and are unlinked in the other, the same with the $2$nd and $4$th. Therefore the Seifert matrix is
$$
\left( \begin{array}{cccc} -n_1 s_1 - n_5 s_5  & n_1 s_1 & 1 & 0\\
n_1 s_1 & -n_1 s_1 - n_4 s_4 & 0 & 1 \\
0 & 0 & -n_2 s_2 & 0 \\
0 & 0 & 0 & -n_3 s_3 \\
\end{array} \right)
$$
Computing the leading term of the Alexander polynomial from this matrix, this is
$$
n_2 n_3 s_2 s_3 (n_1 s_1 n_4 s_4 + (n_4 s_4 + n_1 s_1) n_5 s_5).
$$
The difference between the knots $10_{58}$ and $8_{12}$ is precisely that the $5$th band in $10_{58}$ has zero twists in $8_{12}$ \ref{8-12}. Therefore, the Seifert matrix for the family generated by $8_{12}$ is the same as the one for the family generated by $10_{58}$, except $n_5$ is zero. Making the substitution in the Alexander polynomial, this causes the last two terms to vanish, leaving only a single monomial remaining. Because the $n_i$s are positive integers and the $s_i$s are $\pm 1$ in every case, the leading term of the Alexander polynomial for $8_{12}$ never vanishes, so no possible knot in that family can have cosmetic surgeries.

\section{Computations for the Knot Family Generated by $10_{58}$}
For the knot $10_{58}$, there are again $2^5 = 32$ cases for signs, and all will be handled separately, as in the case of $7_6$, and, again ordered from simplest to longest computation. Unlike in the case of $7_6$, there will be no unknots hidden within the family, so that check will not be required.

The cases are as follows:\\
$(- - - - -)$, $(+ + + + +)$, $(+--++)$, $(-++--)$\\
The leading coefficient of the Alexander polynomial is $ (a d + a e + d e) b c$. This is positive, so nonvanishing. This excludes these cases.\\

$(--+--)$ $(++-++)$, $(+-+++)$, $(-+---)$\\
The leading coefficient of the Alexander polynomial is $-(a d + a e + d e) b c$. This is negative, so nonvanishing. This excludes these cases.\\

$(++++-)$, $(----+)$\\
The leading coefficient of the Alexander polynomial is $(a d - a e - d e) b c$ This vanishes only if $e = \frac{a d}{a + d}$. The second coefficient of the Alexander polynomial is $a b + a c + c d - b e - 4b c(a d - a e - d e)$.  Substituting for $e$ gives
$$
\frac{ (a b + a c + c d)(a + d) - a b d }{a + d} = \frac{a^2 b + (a + d)^2 c}{a + d}
$$
This is clearly always positive, so this case is excluded.\\
Related: $(+--+-)$, $(-++-+)$\\
The difference is that the second coefficient of the Alexander polynomial is $-(a b + a c + c d - b e) - 4b c(a d - a e - d e)$ so substituting for $e$ gives $\frac{-(a^2 b + (a + d)^2 c)}{a + d}$, which is always negative, so cannot vanish and this case is excluded.\\

$(+++-+)$, $(- - - + - )$\\
The leading coefficient of the Alexander polynomial is $(a e - a d - d e) b c$. This vanishes only if $d = \frac{a e}{a + e}$. Considering the second coefficient of the Alexander polynomial and the result of substituting for $d$
$$
a b + a c - c d + b e - 4 (a e - a d - d e) b c = \frac{(a b + a c + b e)(a + e) - c a e}{a + e} = \frac{a^2 c + b(a + e)^2}{a+e}
$$
This is always positive, so this excludes this case.\\
Related: $(+---+)$, $(-+++-)$\\
The difference is that the second coefficient of the Alexander polynomial is $-(ab + ac - c d + b e)  - 4 (a e - a d - d e) b c$ so substituting for $d$ gives $\frac{-(a^2 c + b (a + e)^2)}{a + e}$ which is always negative so this excludes this case.\\

$(+++--)$, $(- - - + +)$\\
The leading coefficient of the Alexander polynomial is $(-a d - a e + d e) b c$. This vanishes only if $a = \frac{d e}{d + e}$. The second coefficient of the Alexander polynomial is $a b + a c - c d - b e - 4 b c(-a d - a e + d e)$. Substituting for $a$ gives $\frac{-(c d^2 + b e^2)}{d + e}$. This is always negative, so this never vanishes.\\
Related: $(+----)$, $(-++++)$\\
The difference is that the second coefficient of the Alexander polynomial is $-(a b + a c - c d - b e) - 4 b c(-a d - a e + d e)$. Substituting for $a$ gives $\frac{(c d^2 + b e^2)}{d + e}$. This is always positive, so this never vanishes.\\

$(++--+)$\\
The leading coefficient of the Alexander polynomial is $(a d - a e + d e) b c$, which vanishes only if $d = \frac{a e}{a + e}$. Considering the second term of the Alexander polynomial and the result of substitution for $d$,
$$
a b - a c + c d + b e - 4(a d - a e + d e) b c = \frac{(a+e)^2 b - c a^2}{a + e}
$$
This vanishes if $c = \frac{b(a + e)^2}{a^2}$. The third derivative of the Jones polynomial is:
$$
18 ((a b - a c + c d + b e) - a^2 b - a b^2 + a^2 c + 2 a b c - a c^2 - 
 2 a c d + c^2 d + c d^2 - 2 a b e - b^2 e - b e^2)
$$
which simplifies to
$$
\frac{-18 b e (2 a^3 + a^2 e + a b e + b e^2)}{a^2}
$$
This is always negative, so never vanishes.\\
Related: $(--++-)$\\
The difference is that the third derivative of the Jones polynomial is:
$$
18 ((a b - a c + c d + b e) -(- a^2 b - a b^2 + a^2 c + 2 a b c - a c^2 - 
 2 a c d + c^2 d + c d^2 - 2 a b e - b^2 e - b e^2))
$$
which simplifies to
$$
\frac{18 b e(2 a^3 + a^2 e + a b e + b e^2)}{a^2}
$$
which is always positive, so never vanishes.\\

$(+-++-)$\\
The leading term of the Alexander polynomial is $(-ad + ae + d e) b c$, which vanishes if $e = \frac{ a d}{a  +d}$. Considering the second term of the Alexander polynomial and the result of substitution for $e$,
$$
-ab + a c + c d + b e - 4 (-ad + ae + d e) b c = \frac{(a + d)^2 c - a^2 b}{a + d}
$$
which vanishes if $b = \frac{c(a + d)^2}{a^2}$. The third derivative of the Jones polynomial simplifies to $\frac{-18 c d (2 a^3 + a^2 d + a c d + c d^2)}{a^2}$ which is always negative, so this excludes this case.\\
Related: $(-+--+)$\\
The difference is that the third derivative of the Jones polynomial simplifies to $\frac{18 c d (2 a^3 + a^2 d + a c d + c d^2)}{a^2}$, which is always positive, so this excludes this case.\\

Here are the cases which also need the fourth derivative, illustrating its use:\\
$(++-+-)$\\
The leading coefficient of the Alexander polynomial is $(-ad + a e + d e) b c$
This vanishes if $e = \frac{ a d }{a+ d}$. The second term of the Alexander polynomial is
$$
a b - a c - c d - b e - 4(-ad + a e + d e) b c = \frac{a^2 b - c(a + d)^2}{a + d}
$$
This vanishes only if $b = \frac{c(a+d)^2}{a^2}$.  The third derivative of the Jones polynomial simplifies to $ \frac{ c d (2a^3 + a^2 d - a c d - c d^2)}{a^2}$, which vanishes only if $2a^3 + a^2 d - a c d - c d^2$ does, or if $a^2(2a + d) = c d(a + d)$, which means that, as neither side can vanish, that $c = \frac{a^2(2a + d)}{d(a+d)}$. Substituting for $c$ $b$ and $e$ in the $4$th derivative evaluated at $1$ gives:
$$
\frac{48 a^3(2 a + d)^2(4a^2 + 4 a d + 5 d^2)}{d^2(a + d)}
$$
This is always positive, so this cannot vanish, and this excludes this case.\\
Related: $(--+-+)$\\
This case differs at the third derivative simplifies to $\frac{-c d (2 a^3 + a^2 d - a c d - c d^2)}{a^2}$, which again is zero only if $c = \frac{a^2(2a + d)}{d(a + d)}$. Substituting for all these variables in the fourth derivative gives the same thing, which again is always positive, so does not vanish.\\

$(+-+-+)$\\
The leading coefficient of the Alexander polynomial is $(ad - a e + d e) b c$. This vanishes if $d = \frac{a e}{a + e}$ The second coefficient of the Alexander polynomial is
$$
-a b + a c - c d - b e - 4 (a d - a e + d e) b c = \frac{a^2 c - b(a + e)^2}{a + e}
$$
This vanishes if $c = \frac{b(a + e)^2}{a^2}$. The third derivative of the Jones polynomial simplifies, after substituting for $d$ and $c$, to $\frac{18 b e (a^2(2 a + e) - (a + e) e b)}{a^2}$, which vanishes only if $b = \frac{a^2(2a + e)}{(a + e)e}$. Substituting for $d, c$ and $b$ in the $4$th derivative of the Jones polynomial gives:
$$
\frac{48 a^3 (2 a + e)^2 (4 a^2 + 4 a e + 5 e^2)}{e^2 (a + e)}
$$
which is always positive.\\
Related: $(-+-+-)$\\
The difference is that the third derivative of the Jones polynomial, after substituting for for $c$ and $d$, simplifies to $\frac{-18 b e (a^2(2 a + e) - (a + e) e b)}{a^2}$, which differs by a total sign and again vanishes only if $b = \frac{a^2(2a + e)}{(a+e)e}$. Substituting for $d, c$ and $b$ in the $4$th derivative gives the same formula as before, which, again, is clearly always positive.\\

$(++---)$\\
The leading coefficient of the Alexander polynomial is $(ad + ae - d e)b c$, which vanishes if $a = \frac{d e}{d + e}$. The second coefficient is:
$$
a b - a c + c d - b e - 4 (a d + a e - d e) b c = \frac{ c d^2 - b e^2}{d + e}
$$
This vanishes if $b = \frac{ c d^2}{e^2}$. Considering the third derivative of the Jones polynomial, this becomes $\frac{18 c d^2 (c (d+e)^2 + (d-e)e^2)}{e^2 (d + e)}$. If this vanishes, then $c = \frac{(e-d) e^2}{(d + e)^2}$, so as $c$ is a positive integer by design, we must have $e > d$. Then substituting for the $4$th derivative of the Jones polynomial, we get
$$
\frac{-48 d^3 (e - d)^2 e^3 (5 d^2 + 6 d e + 5 e^2)}{(d + e)^6}
$$
Therefore, as $e > d$, this cannot vanish and is in fact negative, so this case can be excluded.\\
Related: $(--+++)$\\
The first difference is that the $3$rd derivative of the Jones polynomial is
$$
18((a b - a c + c d - b e) + (a^2 b + a b^2 - a^2 c - 2 a b c + a c^2 + 2 a c d - c^2 d - c d^2 - 2 a b e - b^2 e + b e^2))
$$
However, substituting for $a$ and $b$, this becomes $\frac{-18 c d^2 (c (d+e)^2 + (d-e)e^2)}{e^2 (d + e)}$, differing only by a total sign from the original case, which again vanishes if $c = \frac{e^2 (e - d)}{(d + e)^2}$. As $c$ is a positive number, we again have that $e > d$. Considering the $4$th derivative of the Jones polynomial, with these substitutions, it is the same as before which again must be a negative number, hence it is nonvanishing.\\
Related: $(-+-++)$\\
The first difference between this case and the first is that the second coefficient of the Alexander polynomial is $-(ab - ac + c d - b e) - 4(a d + a e - d e) b c$, which gives the same solution for $b$, that this vanishes precisely when $b = \frac{c d^2}{e^2}$. Considering the $3$rd derivative of the Jones polynomial, this is
$$
18(-(ab - ac + c d - b e) + (-a^2 b + a b^2 + a^2 c - 2 a b c + a c^2 - 2 a c d - c^2 d + c d^2 + 
 2 a b e - b^2 e - b e^2) )
$$
Substituting for $a$ and $b$, this becomes $\frac{-18 c d^2 (e^2 (-d + e) + c (d + e)^2)}{e^2 (d + e)}$. This vanishes precisely when $c = \frac{e^2(d - e)}{(d + e)^2}$ so as $c$ is positive, we have $d > e$. Then, the $4$th derivative of the Jones polynomial becomes:
$$
\frac{-48 d^3 (d - e)^2 e^3 (5 d^2 + 6 d e + 5 e^2)}{(d + e)^6}
$$
which is negative, so this cannot vanish.\\
Related: $(+-+--)$\\
The first difference between this and the case immediately prior is that the $3$rd derivative of the Jones polynomial is
$$
18(-(ab - ac + c d - b e) - (-a^2 b + a b^2 + a^2 c - 2 a b c + a c^2 - 2 a c d - c^2 d + c d^2 + 
 2 a b e - b^2 e - b e^2) )
$$
The simplification differs by a total sign, which gives the same solution for $c$. Substituting for all of this into the $4$th derivative gives the same result as before, which is again non-vanishing and so this case can also be eliminated.

\bibliographystyle{plain}
\bibliography{references}

\end{document}